\documentclass[11pt]{amsart}
\pdfoutput=1
\usepackage{graphicx}
\usepackage{amsmath,amscd}
\usepackage{amsthm}
\usepackage{url}

\newtheorem{thm}{Theorem}
\newtheorem{prop}[thm]{Proposition}

\newtheorem{lem}[thm]{Lemma}

\newtheorem*{main}{Main Theorem}

\theoremstyle{definition}

\theoremstyle{definition}

\theoremstyle{remark}

\begin{document}
\title{The Group of Symmetries of the Tower of Hanoi Graph}
\author{So Eun Park}
\thanks{This work was undertaken at the Columbia University REU program supported by NSF
grant DMS-0739392}
\date{}
\maketitle

\begin{figure}[h]
\begin{center}
    \includegraphics[height=30mm]{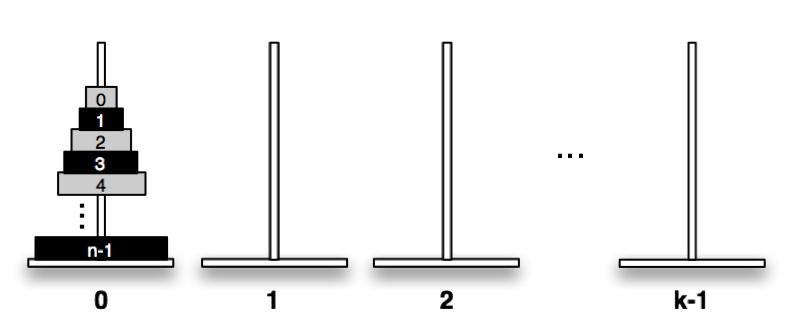}
\end{center}
\caption{Convention for labeling $k$ pegs and $n$ disks in the Tower of Hanoi puzzle.}
  \label{graph labeling}
\end{figure}

The classical Tower of Hanoi puzzle, invented by the French mathematician \'{E}dourd Lucas in 1883, consists of $3$ wooden pegs and $n$ disks with pairwise different diameters. The $n$ disks are initially stacked on a single peg in order of decreasing size, from the largest at the bottom to the smallest at the top. (See Figure~\ref{graph labeling}.) The goal is to move the tower of disks to another peg, moving one topmost disk at a time while never stacking a disk on a smaller one. A sequence of moves realizing this goal in the shortest possible number of moves for any given number of pegs provides a general solution to the Monthly Problem 3918 \cite{Frame}, which is often referred to as the \emph{Tower of Hanoi problem}.

\begin{figure}[t]
\begin{center}
    \includegraphics[height=30mm]{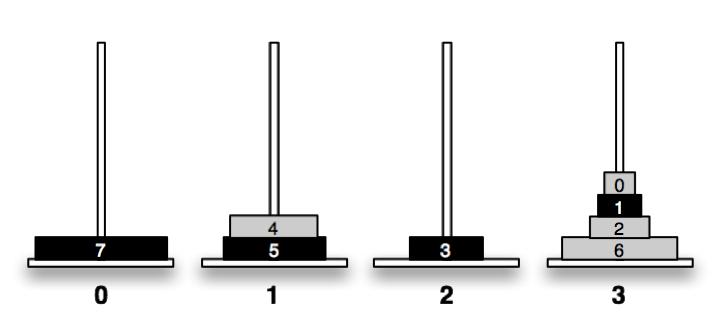}
\end{center}
\caption{An example of a legal state, $a_7a_6a_5a_4a_3a_2a_1a_0=03112333$.}
  \label{graph legal states}
\end{figure}

\begin{figure}[t]
\begin{center}
    \includegraphics[height=30mm]{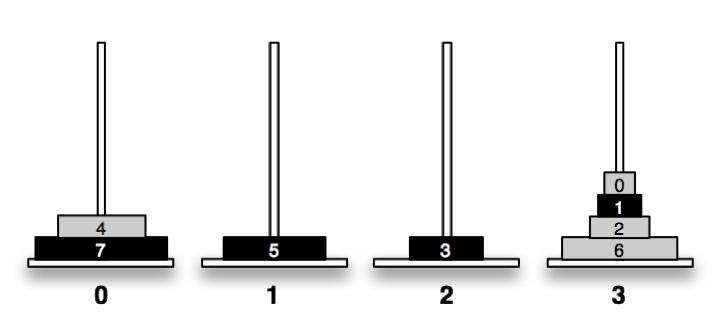}
\end{center}
\caption{A legal state, $a_7a_6a_5a_4a_3a_2a_1a_0=03102333$, that can be reached from Figure~\ref{graph legal states} by moving the topmost disk of peg $1$ to peg $0$.}
  \label{graph legal moves}
\end{figure}

\begin{figure}[t]
\begin{center}
    \includegraphics[height=30mm]{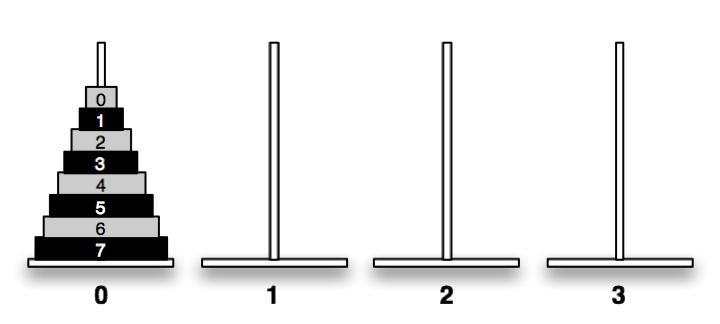}
\end{center}
\caption{An example of a perfect state.}
  \label{graph perfect states}
\end{figure}

\begin{figure}[h]
\begin{center}
    \includegraphics[height=45mm]{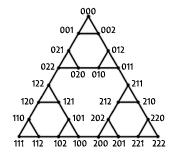}
\end{center}
\caption{The Tower of Hanoi graph, $H_n^k$, for $k=3$ pegs and $n=3$ disks.  Each 3-bit ternary string describes a legal state.}
  \label{graph pascal}
\end{figure}

There have been several discoveries made about the Tower of Hanoi problem and its variations.  The most well-known are the relations between Pascal's triangle and the Tower of Hanoi graph with 3 pegs \cite{Hinz} (See Figure~\ref{graph pascal}.), and an algorithm proposed by Frame and Stewart \cite{Frame} as a solution to the Tower of Hanoi problem, which has not yet been proved to create the shortest path.  The purpose of this note is to introduce a new theorem on the group of symmetries of the Tower of Hanoi graph that may shed further light on solving the Tower of Hanoi problem. In this note, we describe some finite group theory associated to the Tower of Hanoi problem. There have been interesting approaches to the Tower of Hanoi problem coming from geometric group theory, which involve associating certain infinite groups to the game; for details we refer the readers to \cite{Sunik1} and \cite{Sunik2}.

We begin by reviewing some standard definitions for graphs.  Given a {\em graph}, $\Gamma$, the set of vertices is denoted by $V(\Gamma)$, the set of edges by an {\em edge matrix}, $[E(\Gamma)]$ (with an arbitrary but fixed choice of ordering of vertices), where $e_{ij}$ of $[E(\Gamma)]$ is the number of edges between $v_i$ and $v_j$. Two vertices, $v_i$ and $v_j$, are adjacent if $e_{ij}>0$. The {\em degree} of a vertex is defined as $deg(v_i)=\sum_{j=1}^{N} e_{ij}$ where $N=|V(\Gamma)|$.  The distance between vertices is defined to be the length of the shortest edge path between them, i.e., $d(v,v')=\displaystyle\min_{\gamma}\{l(\gamma)\}$ where $\gamma$ ranges over all paths between $v$ and $v'$.  An \emph{automorphism} of $\Gamma$ is an adjacency-preserving (more precisely, edge matrix-preserving) bijection $g: V(\Gamma) \rightarrow V(\Gamma)$, and $G(\Gamma)$ denotes the group of automorphisms of $\Gamma$.  Any automorphism of $\Gamma$ is an isometry of $V(\Gamma)$ with respect to this metric $d$.

In this note, we will borrow most of our notation from \cite{Hinz} and \cite{Poole}, but will use $H_n^k$ to denote the graph associated to the Tower of Hanoi puzzle with $n$ disks and $k$ pegs.  Let us recall that a vertex $v\in V(H_n^k)$ is an $n$-bit $k$-ary string, $a_{n-1}a_{n-2}\cdots a_{0}$, with $a_i\in\{0,\ldots,k-1\}$.  Such a vertex corresponds to the legal state in which disk $i$ lies on peg $a_i$.  (See Figure 2 and Figure 3.)  In particular, for the sake of simplicity we denote by $\overline{i_n}$ the legal state, $ii\cdots i$, in which all disks are stacked on peg $i$.  Such a configuration is called a perfect state.  (See Figure 4.)  Define a \emph{substructure}, $[i]$, to be the set of vertices of $H_n^k$ whose $n$-bit strings correspond to legal states in which the largest disk lies on peg $i$.

$H_n^3$ has a particularly beautiful and simple recursive structure with a fascinating connection to Pascal's triangle, c.f., \cite{Hinz}.  In $H_n^3$, it is straightforward to prove that $G(H_n^3) \cong S_3$, where $S_3$ denotes the permutation group on 3 elements.  This can be done by using the fact that geodesics (shortest paths with respect to the standard graph metric) between perfect states are unique in $H_n^3$.  This implies that any automorphism mapping the perfect states to themselves must also map each perfect state of every recursive substructure to itself.  Therefore, by a simple inductive argument, any automorphism fixing perfect states must be the identity, hence $G(H_n^k)$ is isomorphic to the group of permutations of the perfect states, i.e., $S_3$.

This argument fails when $k > 3$, since it is a well-known result that geodesics are no longer unique.  The purpose of this note is to provide the following extension of the above result to the case $k\geq 3$.

\begin{main}
$G(H_n^k) \cong S_k$ for all $k\geq 3$ and all $n\geq 1$.
\end{main}

We will show that every element of $G(H_n^k)$ is induced by a peg permutation and that the collection of such automorphisms is isomorphic to $S_k$.  For each $\sigma \in S_k$, define the map \[g_{\sigma}: V(H_n^k) \rightarrow V(H_n^k)\] by \[g_{\sigma}(a_{n-1} \cdots a_1a_0) := \sigma(a_{n-1}) \cdots \sigma(a_1)\sigma(a_0)\] for $a_{n-1}\cdots a_1a_0 \in V(H_n^k)$.  Denote by $G(S_k)$ the set $\{g_{\sigma} \,\,|\,\, \sigma \in S_k\}$.  $G(S_k)$ is canonically isomorphic to $S_k$.  In step 1 that follows, we show that any element of $G(S_k)$ is an element of $G(H_n^k)$.  We then show, in step 2, that each element of $G(H_n^k)$ arises as an element of $G(S_k)$, completing the proof of the main theorem.

\subsection*{\bf{Step 1}: $G(S_k) \leq G(H_n^k)$}

\begin{prop} \label{permutation subgroup}
For each $\sigma \in S_k$, $g_\sigma$ is an automorphism of $H_n^k$.  Hence, $G(S_k) \leq G(H_n^k)$.
\end{prop}

\begin{proof}
It is an easy consequence from the physical observation that re-assigning the labels of the pegs does not change the structure of the graph.  Moreover, the automorphism induced by a peg permutation naturally has its inverse that is also induced by another peg permutation.  Therefore, the collection of automorphisms of the Tower of Hanoi graph induced by peg permutations forms a subgroup of $G(H_n^k)$, i.e., $G(S_k) \leq G(H_n^k)$.
\end{proof}

\subsection*{Step 2: $G(H_n^k) \leq G(S_k)$}
$\\$We begin with a useful lemma.

\begin{lem} \label{corner vertices permutation}
Every automorphism in $G(H_n^k)$ permutes the corner vertices, $\{\overline{0_n},\overline{1_n},\ldots,\overline{(k-1)_n}\}$.  That is, for each $g \in G(H_n^k)$ and $i\in\{0,\ldots,k-1\}$, there exists a unique $j\in\{0,\ldots,k-1\}$ such that $g(\overline{i_n}) = (\overline{j_n})$.
\end{lem}

\begin{proof}
We use the fact that the degree of the corner vertices is strictly smaller than that of non-corner vertices.  Since a graph automorphism preserves degree, it must therefore send corner vertices to corner vertices.

To see that the degree of each corner vertex is strictly smaller than the degree of each non-corner vertex, we compute the degree of a vertex in terms of the number of ``topmost disks'' of each vertex.  For a given vertex, the topmost disk on each peg $i$ is defined to be the smallest disk among those stacked on peg $i$.  Therefore, the degree of a vertex is the sum of the number of legal moves each of its topmost disks can make.  Since a corner vertex has only one topmost disk, which is disk $0$, its degree is $k-1$, the number of legal moves disk $0$ can make.  Now we must show that the degree of each non-corner vertex is strictly larger than $k-1$.  If a vertex, $v$, is not a corner vertex, it has $n$ disks distributed on at least $2$ pegs.  Hence, it has at least $2$ topmost disks.  Label the smallest and second smallest of those topmost disks, respectively $b_0$ and $b_1$.  Then $b_0$ can be moved to any other $k-1$ pegs and $b_1$ to any other except where $b_0$ is stacked on, since $b_0$ is the only topmost disk that is smaller than $b_1$.  Therefore, every non-corner vertex has degree of at least $(k-1)+(k-2)=2k-3$, which is strictly larger than $k-1$ for every $k\geq 3$.
\end{proof}

We now prove that the only automorphism in $G(H_n^k)$ that fixes the corner vertices is the identity.

\begin{prop} \label{identity map}
If $g \in G(H_n^k)$ satisfies \[g(\overline{i_n}) = \overline{i_n} \,\,\, \forall\,\, i \in \{0, \ldots, k-1\},\] then $g$ is the identity automorphism.  In other words, the only automorphism in $G(H_n^k)$ that fixes the corner vertices is the identity.
\end{prop}

Proposition~\ref{identity map} implies $G(H_n^k) \leq G(S_k)$ because, by Lemma~\ref{corner vertices permutation}, every $g\in G(H_n^k)$ induces a permutation of the corner vertices, hence of the set $\{0, \ldots, k-1\}$.  Therefore, there exists a $\sigma \in S_k$ that induces the same permutation of the set $\{0, \ldots, k-1\}$ as $g$ does.  Then, $g_{\sigma^{-1}}\circ g(\overline{0_n}, \overline{1_n}, \ldots, \overline{(k-1)_n}) = (\overline{0_n}, \overline{1_n}, \ldots, \overline{(k-1)_n})$.  Hence, Proposition~\ref{identity map} implies that $g_{\sigma^{-1}}\circ g$ is the identity automorphism, making $g=g_{\sigma}\in G(S_k)$.  Hence $G(H_n^k)\leq G(S_k)$.

\begin{proof} [Proof of Proposition~\ref{identity map}]
We proceed by induction on $n$, fixing $k$.
$\\${\bf Base case: $n=1$}.  Trivial, since all vertices of $H_1^k$ are corner vertices and the vertices are connected to one another by a single edge.  Hence, any automorphism which fixes all the corner vertices fixes all vertices as well as edges, thereby inducing the identity automorphism.
$\\${\bf Inductive Step:} Assume that the proposition holds for $n-1$.  To prove that this implies that the proposition holds for $n$, we will need three lemmas, each helping prove the subsequent one.  As we could not find the reference in literature, the argument for the following lemmas is provided by Michael Rand (personal communication, 2008).

\begin{lem} \label{at most once}
Any shortest path between a corner vertex, $\overline{i_n}$, and an arbitrary vertex, $v$, involves moving the largest disk zero times if $v\in[i]$; once otherwise.
\end{lem}

\begin{proof}
By looking at an initial path from a corner vertex, $\overline{i_n}$, and an intermediate vertex, $w$, on a shortest path from  $\overline{i_n}$ to $v$, it suffices to show that there exists no such an initial path that moves the largest disk twice. This suffices since every shortest path between a corner vertex and an arbitrary vertex that moves the largest disk more than zero times, when $v\in[i]$, or more than once, when $v\notin [i]$, has to necessarily contain an initial path that moves the largest disk twice.

There are two cases when the largest disk moves twice from $\overline{i_n}$; the final position of the largest disk is in a different substructure,  $[j]\neq[i]$, or in the same substructure, $[i]$. We will first treat the first case, and then the second case. Without loss of generality, aiming for contradiction, assume that there is a shortest path from $\overline{0_n}$ to $w\in[2]$ moving the largest disk exactly twice---first to peg $1$ and then to peg $2$.  Then we can write the path as a sequence of steps:
\begin{enumerate}
\item Move the $n-1$ smallest disks off of peg $0$ (leaving peg $1$ clear at the end).
\item Move the largest disk from $0$ to $1$.
\item Some number of moves on the $n-1$ smallest disks (maybe 0) which leave peg $2$ empty and peg $1$ containing only the largest disk.
\item Move the largest disk from $1$ to $2$.
\item Some number of moves (maybe 0) to get to the vertex $v$.
\end{enumerate}

We claimed that this is the shortest path from $\overline{0_n}$ to $v$, but we can create an even shorter path as follows.
\begin{enumerate}
\item Do the same moves as in step (1) above, but with the roles of pegs 1 and 2 switched.	
\item Move the largest disk from $0$ to $2$.	
\item Do the same moves as in step (3) above, but with the roles of pegs 1 and 2 switched.	
\item Repeat step $(5)$ above.	
\end{enumerate}

The second sequence of steps is one legal move shorter than the initial sequence, hence contradicts the assumption.  Therefore, there is no such shortest path between a corner vertex and a vertex in a different substructure that involves moving the largest disk exactly twice.

Now, the second case where the largest disk moves twice and comes back to the initial substructure in a shortest path can be treated in a similar fashion. Without loss of generality, assume that there is a shortest path from $\overline{0_n}$ to $w\in[0]$ moving the largest disk exactly twice---first to peg $1$ and then to peg $0$. Then we can write the path as a sequence of steps the way we did above by switching peg $2$ to peg $0$ in step (4). Then, the same sequence of moves without step (2) and (4) gives a new path two legal moves shorter, hence giving a contradiction to our assumption. Therefore, there is no such shortest path between a corner vertex and a vertex in the same substructure that involves moving the largest disk exactly twice.
\end{proof}

\begin{lem} \label{shortest perfect state}
For all $v \in [i],\,\,\, d(v,\overline{i_n}) < d(v,\overline{j_n}) \,\, \forall \,\, j \neq i.$
\end{lem}

\begin{proof}
To simplify notation we re-index the pegs; doing this, we assume for the rest of the argument that $i=0$ and $j=1$.  By Lemma~\ref{at most once}, the shortest path, $\gamma$, from $v\in[0]$ to $\overline{1_n}$ involves moving the largest disk exactly once.  Therefore we can further assume that $\gamma$ can be split into three parts: $\gamma_1$, $\gamma_2$, and $\gamma_3$.  More precisely, $\gamma_1$ is the path from $v$ to $v_0=0a_{n-2}\cdots a_1a_0$, $a_i\notin \{0,1\}$; $\gamma_2$ is the intermediate single legal move that is from $v_0=0a_{n-2}\cdots a_1a_0$ to $v_1=1a_{n-2}\cdots a_1a_0$; finally, $\gamma_3$ is from $v_1=1a_{n-2}\cdots a_1a_0$ to $\overline{1_n}$.  (See Figure~\ref{graph proof}.)  We can easily observe that $l(\gamma_3)=d(v_1,\overline{1_n})=d(v_0,\overline{0_n})$ since the relationship between configurations $v_1$ and $\overline{1_n}$, and that of $v_0$ and $\overline{0_n}$ are exactly symmetric.  Hence, \begin{eqnarray*}
  d(v,\overline{0_n}) &\leq& d(v,v_0) + d(v_0, \overline{0_n}) \,\,\,\,\,\,\,\,\, \mbox{ (by the triangle inequality)}\\
           &<& d(v,v_0) + 1 + d(v_0,\overline{0_n})\\
           &=& d(v,v_0) + 1 +d(v_1,\overline{1_n})\\
           &=& l(\gamma_1) + l(\gamma_2) + l(\gamma_3)\\
           &=& l(\gamma)\\
           &=& d(v,\overline{1_n}).
\end{eqnarray*}
\end{proof}

\begin{figure}
\begin{center}
    \includegraphics[height=45mm]{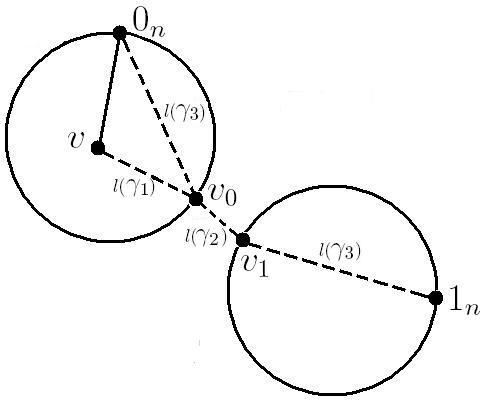}
\end{center}
\caption{The sketch of proof for Lemma~\ref{shortest perfect state}. $l(\gamma_2)=d(v_1,\overline{1_n})=d(v_0,\overline{0_n})$ and triangle inequality is applied to find $d(v,\overline{0_n})<d(v,\overline{1_n})$.}
  \label{graph proof}
\end{figure}

\begin{lem} \label{substructure preservation}
If $g \in G(H_n^k)$ satisfies $g(\overline{i_n}) = \overline{i_n}$, then $g(v) \in [i]$ for all $v \in [i]$.  I.e., $g([i]) = [i]$ as a set.
\end{lem}

\begin{proof}
For concreteness we assume that $\overline{0_n}$ is the vertex which is fixed by hypothesis.  Now assume, aiming for a contradiction, that there exists $v \in [0]$ such that $g(v) \in [i], i \neq 0$.  By Lemma~\ref{corner vertices permutation}, there exists a unique $j \neq 0$ such that $g(\overline{j_n}) = \overline{i_n}$.  Thus, both $g(v)$ and $g(\overline{j_n})$ are in $[i]$. Thus, since
$g(\overline{j_n})=\overline{i_n}$, by Lemma~\ref{shortest perfect state} we have \[d(g(v),g(\overline{j_n})) < d(g(v), g(\overline{k_n}))\] $\forall k \neq j$. Since $\forall g\in G(H_n^k)$ is an isometry, the previous statement implies $d(v,\overline{j_n})<d(v,\overline{k_n}),\,\,\forall k\neq j$.  Since $0\neq j$, it follows that $d(v,\overline{j_n})<d(v,\overline{0_n})$ with $v\in [0]$, which contradicts Lemma~\ref{shortest perfect state}.
\end{proof}

Now assume that Proposition~\ref{identity map} holds for $G(H_{n-1}^{k})$ for all $k\geq 3$, and let $g\in G(H_n^k)$ satisfy $g(\overline{i_n}) =\overline{i_n}$ for all $i \in \{0, \ldots, k-1\}$.  We have to show that this implies that $g$ is the identity automorphism of $H_n^k$.  By Lemma~\ref{substructure preservation}, $g(\overline{i_n}) =\overline{i_n}$ for all $i \in \{0, \ldots, k-1\}$ implies $g([i])=[i]$ for all $i \in \{0, \ldots, k-1\}$.  Thus, for all $i$, $g|_{[i]}$ is an automorphism of \[\{a_{n-1}a_{n-2}\cdots a_1a_0 \in V(H_{n-1}^k)\,\,|\,\,a_{n-1}=i\}\]  Since the leading entry, $i$, is fixed, the automorphism $g$ restricted to $[i]$ induces an automorphism, $g_{i}:V(H^{k}_{n-1})\to V(H^{k}_{n-1})$, satisfying $g_{i}(\overline{i_{n-1}})=\overline{i_{n-1}}$.

By Lemma~\ref{corner vertices permutation}, for each $j\neq i$, there exists a unique $l\neq i$ such that $g_i(\overline{j_{n-1}})=\overline{l_{n-1}}$, hence $g(i\overline{j_{n-1}})=i\overline{l_{n-1}}$.  We will now show that $j=l$, hence $g(i\overline{j_{n-1}})=i\overline{j_{n-1}}$ for all $j$.  We note the observation that $i\overline{j_{n-1}}$ is never adjacent to any vertex in $[j]$, but on the other hand, has an adjacent vertex, $l\overline{j_{n-1}}$, in any other substructure $[l]$, $\forall l\neq j$.  Since $g$ is an automorphism it preserves adjacency, hence the above observation implies that $g(i\overline{j_{n-1}})$ is still never adjacent to any vertex in $g([j])=[j]$, but has an adjacent vertex in any other substructure $g([l])=[l]$, $\forall l\neq j$.  Hence $g(i\overline{j_{n-1}})=i\overline{j_{n-1}}$ since otherwise it must be adjacent to a vertex in $[j]$. Thus, $g(i\overline{j_{n-1}})=i\overline{j_{n-1}}$ for all $j$.  This implies that for all $i$, $g_i \in G(H_{n-1}^k)$ satisfies the assumption of Proposition~\ref{identity map}.  Hence, by the inductive hypothesis, $g_i$ is the identity for all $i$, implying that $g$ fixes all the vertices of $H^k_n$.  Thus, since there is at most one edge between any pair of vertices, any automorphism which fixes all the vertices must fix all the edges as well.  Thus, we have shown $g$ is the identity automorphism.
\end{proof}

\textbf{Step 1} and \textbf{Step 2} are now proved.  Hence, $G(H_n^k) \cong G(S_k) \cong S_k$, as desired.

\paragraph{Acknowledgments.} I would like to thank Professors Jason Behrstock and Elisenda Grigsby, the graduate student advisor, Harold Sultan, and the participants in the Columbia University 2008 summer REU program, especially Michael Rand for offering the figures for the proof as well as the idea for the proof of Lemma 4.  I would like to thank Columbia University math department, Professor Robert Friedman and Peter Ozsv\'{a}th, and the National Science Foundation for organizing and supporting the program.  I would also like to thank the referees for their useful comments about revision.

\bigskip

\noindent\textit{sp2394@columbia.edu}

\end{document}